\documentclass[11pt]{article}
\usepackage{amsfonts}
\usepackage{latexsym, amssymb, amsmath, amscd, amsfonts, epsfig, graphicx, colordvi}
\usepackage{ifpdf}
\usepackage{graphicx}
\usepackage{xcolor}
\usepackage[dvipsnames]{pstricks}
\newtheorem{thm}{Theorem}[section]

\newtheorem{lem}[thm]{Lemma}

\newtheorem{core}[thm]{Corollary}

\def\pf{\noindent{\it Proof.} }
\def\qed{\nopagebreak\hf^{}ill{\rule{4pt}{7pt}}
\medbreak}

\def\qed{\nopagebreak\hfill{\rule{4pt}{7pt}}
\medbreak}

\makeatletter
\def\ExtendSymbol#1#2#3#4#5{\ext@arrow 0099{\arrowfill@#1#2#3}{#4}{#5}}
\makeatother

\begin{document}

\begin{center}
{\large\bf An Inequality for Coefficients of the Real-rooted Polynomials}
\end{center}
\begin{center}
	Jeremy J. F. Guo  \\[8pt]
	School of Mathematics\\
	Tianjin University\\
	Tianjin 300072, P. R. China\\[6pt]

	Email: {\tt guo@tju.edu.cn}, 
\end{center}

\vspace{0.3cm} \noindent{\bf Abstract}
In this paper, we prove that if $f(x)=\sum_{k=0}^n{n\choose k}a_kx^k$ is a polynomial with real zeros only, then the sequence $\{a_k\}_{k=0}^n$ satisfies the following inequalities $a_{k+1}^2(1-\sqrt{1-c_k})^2/a_k^2
\leq(a_{k+1}^2-a_ka_{k+2})/(a_k^2-a_{k-1}a_{k+1})
\leq a_{k+1}^2(1+\sqrt{1-c_k})^2/a_k^2$,
where $c_k=a_ka_{k+2}/a_{k+1}^2$. This inequality holds for the coefficients of the Riemann $\xi$-function, the ultraspherical, Laguerre and Hermite polynomials, and the partition function. Moreover, as a corollary, for the partition function $p(n)$, we prove that $p(n)^2-p(n-1)p(n+1)$ is increasing for $n\geq 55$. We also find that for a positive and log-concave sequence $\{a_k\}_{k\geq 0}$, the inequality $a_{k+2}/a_k\leq (a_{k+1}^2-a_ka_{k+2})/(a_k^2-a_{k-1}a_{k+1}) \leq a_{k+1}/a_{k-1}$ is the sufficient condition for both the  $2$-log-concavity and the higher order Tur{\'a}n inequalities of $\{a_k\}_{k\geq 0}$. It is easy to verify that if $a_k^2\geq ra_{k+1}a_{k-1}$, where $r\geq 2$, then the sequence $\{a_k\}_{k\geq 0}$ satisfies this inequality.



\section{Introduction}

In this paper, we give an inequality for the coefficients of real-rooted polynomials. 

\begin{thm}\label{mainthm}
For a real-rooted polynomial $f(x)=\sum_{k=0}^n {n\choose k}a_kx^k$, if  $a_ka_{k+1}(a_k^2-a_{k-1}a_{k+1})\neq0$, then the inequality
\begin{equation}\label{1}
\frac{a_{k+1}^2}{a_k^2}(1-\sqrt{1-c_k})^2
\leq\frac{a_{k+1}^2-a_ka_{k+2}}{a_k^2-a_{k-1}a_{k+1}}
\leq
\frac{a_{k+1}^2}{a_k^2}(1+\sqrt{1-c_k})^2
\end{equation}
holds for $1\leq k\leq n-2$, where $c_k=\frac{a_ka_{k+2}}{a_{k+1}^2}$.
\end{thm}

 We say a polynomial $f(x)=\sum_{k=0}^na_kx^k$ is {\it real-rooted}, if all its zeros are real. The inequality (\ref{1}) gives an upper and lower bound for the ratio $\frac{a_{k+1}^2-a_ka_{k+2}}{a_k^2-a_{k-1}
a_{k+1}}$. And it is equivalent to the higher order Tur{\'a}n inequality of the sequence  $\{a_k\}_{k\geq 0}$. 


Recall that a sequence $\{a_k\}_{k\geq 0}$ is said to be {\it log-concave}  if for all $k\geq 1$,
\begin{equation}\label{lcav}
a_k^2-a_{k-1}a_{k+1}\geq 0.
\end{equation}
Note that for a positive sequence $\{a_k\}_{k\geq 0}$, it is log-concave if and only if the ratio $a_{k+1}/a_k$ is decreasing. We also say that the sequence $\{a_k\}_{k\geq 0}$ satisfies the {\it Tur{\'a}n inequalities}, if it satisfies the inequality (\ref{lcav}).

 For the sequence $\{a_k\}_{k\geq 0}$ satisfying the Tur{\'a}n inequalities, we consider the higher order Tur{\'a}n inequalities as follows. A sequence $\{a_k\}_{k\geq 0}$ is said to satisfy the {\it higher order Tur{\'a}n inequalities} if for $k\geq 1$,
\begin{equation}\label{ht}4(a_k^2-a_{k-1}a_{k+1})(a_{k+1}^2-a_ka_{k+2})-
(a_ka_{k+1}-a_{k-1}a_{k+2})^2\geq 0.
\end{equation}

Recall that a real entire function
\begin{equation}
\psi(x)=\sum_{k=0}^{\infty}\gamma_k
\frac{x^k}{k!}
\end{equation}
is said to be in the Laguerre-P{\'o}lya class, denoted $\psi(x)\in \mathcal{LP}$, if it can be represented in the form
\begin{equation}
\psi(x)=cx^me^{-\alpha x^2+\beta x}\prod_{k=1}^{\infty}(1+x/x_k)e^{-x/x_k},
\end{equation}
where $c$, $\beta$, $x_k$ are real numbers, $\alpha\geq 0$, $m$ is a nonnegative integer and $\sum x_k^{-2}<\infty$. These functions are only ones which are uniform limits of polynomials whose zeros are real. We refer to \cite{lev} and \cite{rah} for the background on the theory of the $\mathcal{LP}$ class.

For a real entire function $\psi(x)=\sum_{k=0}^{\infty}\gamma_k
\frac{x^k}{k!}$ in the $\mathcal{LP}$ class, the Maclaurin coefficients $\gamma_k$ satisfy both the Tur{\'a}n inequalities, proved by P{\'o}lya and Schur \cite{pol}, and the higher order Tur{\'a}n inequalities, proved by  Dimitrov \cite{dim1}. As a corollary, the ultraspherical, Laguerre and Hermite polynomials satisfy both the Tur{\'a}n inequalities and the higher order Tur{\'a}n inequalities, see \cite{dim1}.

Since the inequality (\ref{1}) is equivalent to the higher order Tur{\'a}n inequality, then we get that for a real entire function $\psi(x)$ in the $\mathcal{LP}$ class, the Maclaurin coefficients satisfy the inequality (\ref{1}). Consequently, the ultraspherical, Laguerre and Hermite polynomials satisfy the inequality (\ref{1}).

To prove the higher order Tur{\'a}n inequalities for the Maclaurin coefficients, Dimitrov applied a theorem of Ma{\u r}{\'i}k \cite{mar} as follows.
\begin{thm}\label{marik}
	If the real polynomial
	$f(x)=\sum_{k=0}^n a_kx^k/(k!(n-k)!)$
	of degree $n\geq 3$ has only real zeros, then the inequality
	\begin{equation*}4(a_k^2-a_{k-1}a_{k+1})(a_{k+1}^2-a_ka_{k+2})-
	(a_ka_{k+1}-a_{k-1}a_{k+2})^2\geq 0
	\end{equation*}
	holds for $1\leq k\leq n-1$.
\end{thm}

It is well known that the Riemann hypothesis holds if and only if the Riemann $\xi$-function belongs to the $\mathcal{LP}$ class.
Let $\zeta$ denote the Riemann zeta-function and $\Gamma$ be the gamma-function. The Riemman $\xi$-function is defined by
\begin{equation}
\xi(iz)=\frac{1}{2}(z^2-\frac{1}{2})
\pi^{-z/2-1/4}\Gamma(\frac{z}{2}+
\frac{1}{4})\zeta(z+\frac{1}{2}),
\end{equation}
see, for example, Boas \cite{boa}.  Hence, if the Riemann hypothesis is true, then the Maclaurin coefficients of the Riemann $\xi$-function satisfy both the Tur{\'a}n inequalities and the higher order Tur{\'a}n inequalities. Csordas, Norfolk and Varga \cite{ccv} proved that the coefficients of the Riemann $\xi$-function satisfy the Tur{\'a}n inequalities, confirming a conjecture of P{\'o}ly \cite{pol1}. Dimitrov and Lucas \cite{dim2} showed that the coefficients of the Riemann $\xi$-function satisfy the higher order Tur{\'a}n inequalities without resorting to the Riemann hypothesis. As a corollary, we conclude that the coefficients of the Riemann $\xi$-function satisfy the inequality (\ref{1}).

For the partition function $p(n)$, Chen, Jia and Wang \cite{cjw} proved that it satisfies the higher order Tur{\'a}n inequalities for $n\geq 95$. As a corollary, the inequality (\ref{1}) holds for partition function $p(n)$ for $n\geq 95$.
 
Through the disscussion  about the lower bound $a_{k+1}^2(1-\sqrt{1-c_k})^2/a_k^2$ in the inequality (\ref{1}), we prove that for the partition function  $p(n)$, $p(n)^2-p(n-1)p(n+1)$ is increasing for $n\geq 55$.
 
Go back to the log-concavity of the sequence $\{a_k\}_{k\geq 0}$. We consider the {\it$2$-log-concavity}, which is equivalent to the decreasing property of the ratio $\frac{a_{k+1}^2-a_ka_{k+2}}{a_k^2-a_{k-1}
	a_{k+1}}$. Moreover, we could define the {\it infinitely  log-concave sequence } as follows.

Define an operator $\mathcal{L}$ on a sequence $\{a_k\}_{k\geq 0}$ by
$\mathcal{L}(\{a_k\}_{k\geq 0})=\{b_k\}_{k\geq 0}$,
where $b_0=a_0^2$ and $b_k=a_k^2-a_{k-1}a_{k+1}$. This definition makes sense for finite sequences by regarding these as infinite sequences
with finitely many nonzero entries. Hence a sequence $\{a_k\}_{k\geq 0}$ is log-concave if and only if $\mathcal{L}(\{a_k\}_{k\geq 0})$ is a nonnegative sequence. We say that a sequence $\{a_k\}_{k\geq 0}$ is {\it $k$-log-concave} if $\mathcal{L}^j(\{a_k\}_{k\geq 0})$ is nonnegative for all $0\leq j\leq k$. A sequence $\{a_k\}_{k\geq 0}$ is {\it infinitely log-concave} if it is $k$-log-concave for all $k\geq 1$.

The notion of infinite log-concavity was introduced by Boros and Moll \cite{bor}. For the sequence $\{{n\choose k}\}_{k=0}^n$, they asked whether it is infinitely log-concave. The following result was independently conjectured by Fisk \cite{fis}, McNamara-Sagan \cite{mcn} and Stanley \cite{sta}, and proved by Br{\"a}nd{\'e}n \cite{bra}.

\begin{thm}\label{branden}
	If $f(x)=\sum_{k=0}^n a_kx^k$ is a polynomial with real- and nonpositive zeros only, then so is
	$$\mathcal{L}(f)=\sum_{k=0}^n
	(a_k^2-a_{k+1}a_{k-1})x^k.$$
	In particular, the sequence $\{a_k\}_{k=0}^n$ is infinitely log-concave.
\end{thm}

It follows immediately that the sequence $\{{n\choose k}\}_{k=0}^n$ is infinitely log-concave.

There is also a simple criterion \cite{cra,mcn} that if
\[
a_k^2\geq ra_{k-1}a_{k+1},\ \  for\ all\ k\geq 1,
\]
where $r\geq (3+\sqrt{5})/2\approx 2.62$, then the sequence $\{a_k\}_{k\geq0}$ is infinitely log-concave.

We are interested in the connection between the $2$-log-concavity and the higher order Tur{\'a}n inequalities. Based on the inequality (\ref{1}), if we can find sharper bounds $l(n)$ and $u(n)$ for the ratio $\frac{a_{k+1}^2-a_ka_{k+2}}{a_k^2-a_{k-1}a_{k+1}}$ such that
\[
\frac{a_{k+1}^2}{a_k^2}(1-\sqrt{1-c_k})^2\leq l(k)
\leq\frac{a_{k+1}^2-a_ka_{k+2}}{a_k^2-a_{k-1}a_{k+1}}
\leq u(k) \leq
\frac{a_{k+1}^2}{a_k^2}(1+\sqrt{1-c_k})^2,
\]
and for each $k\geq 1$, either $\frac{a_{k+1}^2-a_ka_{k+2}}{a_k^2-a_{k-1}a_{k+1}}\geq u(k+1)$
or $\frac{a_{k+2}^2-a_{k+1}a_{k+3}}{a_{k+1}^2-a_ka_{k+2}}\leq l(k)$, then the sequence $\{a_k\}_{k\geq 0}$ is $2$-log-concave, as well as satisfies the higher order Tur{\'a}n inequalities.

In Section 3, we will prove the following theorem.
\begin{thm}\label{2l}
  For a log-concave positive sequence $\{a_k\}_{k\geq0}$, if it satisfies the following inequalities
  \begin{equation}\label{2}
    \frac{a_{k+2}}{a_k}
    \leq\frac{a_{k+1}^2-a_ka_{k+2}}{a_k^2-a_{k-1}a_{k+1}}
    \leq \frac{a_{k+1}}{a_{k-1}},
  \end{equation}
  for $k\geq 1$. Then $\{a_k\}_{k\geq0}$ is $2$-log-concave and satisfies the higher order Tur{\'a}n inequalities for $k\geq 1$.
\end{thm}
It is easy to verify that the sequence $\{{n\choose k}\}_{k=0}^n$ satisfies inequality (\ref{2}), as well as the sequence $\{a_k\}_{k\geq0}$, which satisfies $a_k^2\geq ra_{k+1}a_{k-1}$, where $r\geq 2$.

Finally, in Section 4, we will discuss a problem we will consider in the further work.

\section{Main Theorem}

%
In this section, we will give the proof of the Theorem \ref{mainthm}.

\pf Applying Theorem \ref{marik}, we have
\begin{equation}4(a_k^2-a_{k-1}a_{k+1})(a_{k+1}^2-a_ka_{k+2})\geq
(a_ka_{k+1}-a_{k-1}a_{k+2})^2.
\end{equation}
Multiplying both sides of $1/a_k^2a_{k+1}^2$ and simplifying, we obtain
\begin{equation}4(1-\frac{a_{k-1}a_{k+1}}{a_k^2})(1-\frac{a_ka_{k+2}}{a_{k+1}^2})
\geq(1-\frac{a_{k-1}a_{k+1}}{a_k^2}\frac{a_ka_{k+2}}{a_{k+1}^2})^2.
\end{equation}
Substitute $c_n=\frac{a_na_{n+2}}{a_{n+1}^2}$ for $n=k$ and $n=k-1$, and we get
\begin{equation}
4(1-c_{k-1})(1-c_k)\geq (1-c_{k-1}c_k)^2.
\end{equation}
Observe that
\begin{equation}
1-c_{k-1}c_k=1-c_k+c_k-c_{k-1}c_k
=1-c_k+c_k(1-c_{k-1}).
\end{equation}
It follows that
\begin{equation}4(1-c_{k-1})(1-c_k)\geq
(1-c_k+c_k(1-c_{k-1}))^2.
\end{equation}
Since $a_k^2-a_{k-1}a_{k+1}\neq0$, $1-c_{k-1}\neq0$. Multiply both sides
of $1/(1-c_{k-1})^2$, and we get
\begin{equation}4\frac{1-c_k}{1-c_{k-1}}\geq (\frac{1-c_k}{1-c_{k-1}}+c_k)^2.
\end{equation}
Set $x=(1-c_k)/(1-c_{k-1})$. Then we obtain
\begin{equation}
x^2-(4-2c_k)x+c_k^2\leq 0.
\end{equation}
Immediately we conclude that
\begin{equation}
\frac{4-2c_k-\sqrt{(4-2c_k)^2-4c_k^2}}{2}\leq x
\leq \frac{4-2c_k+\sqrt{(4-2c_k)^2-4c_k^2}}{2}.
\end{equation}
It leads to the following inequality after simplified
\begin{equation}
(1-\sqrt{1-c_k})^2\leq x\leq (1+\sqrt{1-c_k})^2.
\end{equation}
Since
\begin{equation}
\frac{a_{k+1}^2-a_ka_{k+2}}{a_k^2-a_{k-1}a_{k+1}}
=\frac{a_{k+1}^2}{a_k^2}\frac{1-c_k}{1-c_{k-1}}
=\frac{a_{k+1}^2}{a_k^2}\cdot x,
\end{equation}
we reach the inequality (\ref{1}) as we want.\qed

Through the proof above, we can easily get that if $a_ka_{k+1}(a_k^2-a_{k-1}a_{k+1})\neq0$, the inequality (\ref{1}) is equivalent to the higher order Tur{\'a}n inequality. As corollaries, we get the following results.

\begin{core}
  The inequality (\ref{1}) holds for the ultraspherical, Laguerre and Hermite polynomials
\end{core}

\begin{core}
  The inequality (\ref{1}) holds for the coefficients of the Riemann $\xi$-function.
\end{core}

\begin{core}\label{pn}
  The inequality (\ref{1}) holds for the partition function $p(n)$ for $n\geq 95$.
\end{core}

Recall that a sequence $\{a_k\}_{k\geq 0}$ is said to be {\it convex} if for $k\geq 1$,
\begin{equation}
a_{k+1}-a_k\geq a_k-a_{k-1}.
\end{equation}
For the lower bound function $l(n)=\frac{a_{n+1}^2}{a_n^2}(1-\sqrt{1-c_n})^2$ in the inequality (\ref{1}), we have the following result.

\begin{lem}\label{lem21}
	For the log-concave, increasing, positive sequence $\{a_k\}_{k\geq0}$ which satisfies the inequality (\ref{1}), if $\{a_k\}_{k\geq0}$ is convex, then
	\begin{equation}
	\frac{a_{k+1}^2}{a_k^2}(1-\sqrt{1-c_k})^2\geq 1.
	\end{equation}
\end{lem} 
\pf Since
\[
\frac{a_{k+1}^2}{a_k^2}(1-\sqrt{1-c_k})^2
=\frac{1}{a_k^2}(a_{k+1}-a_{k+1}\sqrt{1-c_k})^2
=\frac{1}{a_k^2}(a_{k+1}-\sqrt{a_{k+1}^2-a_ka_{k+2}})^2,
\]
we only need to prove that
\begin{equation}
(a_{k+1}-\sqrt{a_{k+1}^2-a_ka_{k+2}})^2\geq a_k^2.
\end{equation}
For $\{a_k\}_{k\geq0}$ is an increasing, positive sequence, it is sufficient to prove that
\begin{equation}\label{221}
a_{k+1}^2-a_ka_{k+2}\leq (a_{k+1}-a_k)^2.
\end{equation}
Since $\{a_k\}_{k\geq0}$ is convex, for $k\geq 0$
\begin{equation}
a_{k+2}-a_{k+1}\geq a_{k+1}-a_k.
\end{equation}
Thus
\begin{equation}
a_{k+2}\geq 2a_{k+1}-a_k.
\end{equation}
It follows that
\begin{equation}
a_ka_{k+2}\geq a_k(2a_{k+1}-a_k).
\end{equation}
Since
\begin{equation}
a_{k+1}^2-a_k(2a_{k+1}-a_k)=(a_{k+1}-a_k)^2,
\end{equation}
immediately we get the inequality (\ref{221}).
\qed

Combining the Theorem \ref{mainthm} and the Lemma \ref{lem21}, we get the following theorem.
\begin{thm}\label{thm2}
	For the log-concave, increasing, positive sequence $\{a_k\}_{k\geq0}$ which satisfies the inequality (\ref{1}), if $\{a_k\}_{k\geq0}$ is convex, then the sequence $\{a_{k+1}^2-a_ka_{k+2}\}_{k\geq0}$ is increasing.
\end{thm}

For the partition function $p(n)$, $p(n)$ satisfies the inequality \cite{hon}
\begin{equation}
2p(n)\leq p(n+1)+p(n-1).
\end{equation}
Hence we have the corollary as follows.
\begin{core}
	For the partition function $p(n)$, $p(n)^2-p(n-1)p(n+1)$ is increasing for $n\geq 55$.
\end{core}
\pf Applying the Theorem \ref{thm2} and the Corollary \ref{pn}, we get that for $n\geq 95$, $p(n)^2-p(n-1)p(n+1)$ is increasing. For $55\leq n\leq 95$, we can easily verify that $p(n)^2-p(n-1)p(n+1)$ is also increasing.\qed

\section{$2$-log-concavity}

 In this section, we will give the proof of the Theorem \ref{2l}.


\pf 
 Applying the inequalities (\ref{2}), it follows immediately that 
 \begin{equation}
 \frac{a_{k+1}^2-a_{k}a_{k+2}}{a_k^2-a_{k-1}a_{k+1}}
\leq \frac{a_{k+1}}{a_{k-1}}
\leq \frac{a_{k}^2-a_{k-1}a_{k+1}}{a_{k-1}^2-a_{k-2}a_{k}}
\end{equation}
for $k\geq 1$.
 Hence $\{a_k\}_{k\geq0}$ is $2$-log-concave.

On the other hand, for $k\geq 1$, consider the right inequality of (\ref{2}) and we have
\begin{equation}\label{21}
  0\leq\frac{a_{k+1}^2-a_ka_{k+2}}{a_k^2-a_{k-1}a_{k+1}}
  \leq \frac{a_{k+1}}{a_{k-1}}.
\end{equation}
Hence
\begin{equation}
a_{k-1}(a_{k+1}^2-a_ka_{k+2})\leq a_{k+1}(a_k^2-a_{k-1}a_{k+1}).
\end{equation}
Since $a_{k+1}>0$, multiply both sides of $a_{k+1}$, and  we obtain
\begin{equation}\label{22}
a_{k-1}a_{k+1}(a_{k+1}^2-a_ka_{k+2})\leq a_{k+1}^2(a_k^2-a_{k-1}a_{k+1}).
\end{equation}
It leads to
\begin{equation}
-a_{k-1}a_{k}a_{k+1}a_{k+2}\leq a_k^2a_{k+1}^2-2a_{k+1}^3a_{k-1}.
\end{equation}
Thus
\[
-a_{k-1}a_{k}a_{k+1}a_{k+2}+a_k^2a_{k+1}^2\leq a_k^2a_{k+1}^2-2a_{k+1}^3a_{k-1}+a_k^2a_{k+1}^2,
\]
i.e.
\begin{equation}\label{24}
a_ka_{k+1}(a_ka_{k+1}-a_{k-1}a_{k+2})\leq 2a_{k+1}^2(a_k^2-a_{k-1}a_{k+1}).
\end{equation}
Similarly, for $k\geq 1$, consider the left inequality of (\ref{2}) and we get the following inequality
\begin{equation}\label{25}
a_ka_{k+1}(a_ka_{k+1}-a_{k-1}a_{k+2})\leq 2a_{k}^2(a_{k+1}^2-a_ka_{k+2}).
\end{equation} 
Note that $\{a_k\}_{k\geq 0}$ is log-concave. It is easy to verify that
\begin{equation}
a_ka_{k+1}-a_{k-1}a_{k+2}\geq 0.
\end{equation}
Consequently, combine inequalities (\ref{24}) and (\ref{25}), and we get 
\begin{equation}
a_k^2a_{k+1}^2(a_ka_{k+1}-a_{k-1}a_{k+2})^2\leq 4a_k^2a_{k+1}^2(a_k^2-a_{k-1}a_{k+1})(a_{k+1}^2-a_ka_{k+2}).
\end{equation}
 Hence $\{a_k\}_{k\geq 0}$
satisfies the higher order Tur{\'a}n  inequalities.\qed

For a log-concave positive sequence $\{a_k\}_{k\geq 0}$, to prove the inequality (\ref{2}), we need the following two lemmas. 

\begin{lem}\label{lem1}
For a log-concave positive sequence $\{a_k\}_{k\geq 0}$, if the sequence
$\{\frac{a_{k+1}}{a_k}\}_{k\geq 0}$	is convex, then 
\[
\frac{a_{k+1}^2-a_{k}a_{k+2}}{a_k^2-a_{k-1}a_{k+1}}
\leq \frac{a_{k+1}}{a_{k-1}}
\]
for $k\geq 1$.
\end{lem}
\pf Since $\{a_k\}_{k\geq 0}$ is log-concave, $\{\frac{a_{k+1}}{a_k}\}_{k\geq 0}$	is decreasing. Then the convexity of $\{\frac{a_{k+1}}{a_k}\}_{k\geq 0}$ leads to
\[
0\geq\frac{a_{k+2}}{a_{k+1}}-\frac{a_{k+1}}{a_k}\geq\frac{a_{k+1}}{a_k}-\frac{a_k}{a_{k-1}},
\]
i.e.
\begin{equation}
0\leq\frac{a_{k+1}}{a_k}-\frac{a_{k+2}}{a_{k+1}}\leq\frac{a_k}{a_{k-1}}-\frac{a_{k+1}}{a_k}.
\end{equation}
Observe that
\begin{equation}
a_{k+1}^2-a_{k}a_{k+2}=a_{k+1}a_{k}(\frac{a_{k+1}}{a_k}-\frac{a_{k+2}}{a_{k+1}}),
\end{equation}
and
\begin{equation}
a_k^2-a_{k-1}a_{k+1}=a_ka_{k-1}(\frac{a_k}{a_{k-1}}-\frac{a_{k+1}}{a_k}).
\end{equation}
It follows immediately that
\[
\frac{a_{k+1}^2-a_{k}a_{k+2}}{a_k^2-a_{k-1}a_{k+1}}
\leq \frac{a_{k+1}}{a_{k-1}}.
\]
\qed

Similarly, we have the following lemma.
\begin{lem}
	For a log-concave positive sequence $\{a_k\}_{k\geq 0}$, if the sequence
	$\{\frac{a_k}{a_{k+1}}\}_{k\geq 0}$	is convex, then 
	\[
	\frac{a_{k+1}^2-a_{k}a_{k+2}}{a_k^2-a_{k-1}a_{k+1}}
	\geq \frac{a_{k+2}}{a_k}
	\]
	for $k\geq 1$.
\end{lem}

Now we are ready to prove that the sequence $\{{n\choose k}\}_{k=0}^n$ satisfies the inequality (\ref{2}).
\begin{thm}
The sequence	$\{{n\choose k}\}_{k=0}^n$ satisfies the inequality (\ref{2}).
\end{thm}

\pf  It's easy to prove that
\begin{equation}
{n\choose k+1}=\frac{n-k}{k+1}{n\choose k}.
\end{equation}
Since
\begin{equation}
\frac{n-k}{k+1}-\frac{n-k+1}{k}=-\frac{n+1}{k(k+1)}\geq 
-\frac{n+1}{(k-1)k}=\frac{n-k+1}{k}-\frac{n-k+2}{k-1},
\end{equation}
and
\begin{equation}
\frac{k+1}{n-k}-\frac{k}{n-k+1}=\frac{n+1}{(n-k)(n-k+1)}\geq 
\frac{n+1}{(n-k+1)(n-k+2)}=\frac{k}{n-k+1}-\frac{k-1}{n-k+2},
\end{equation}
we conclude that the sequences $\{\frac{n-k}{k+1}\}_{k\geq0}^n$ and $\{\frac{k+1}{n-k}\}_{k\geq0}^n$ are both convex.
Hence the sequence $\{{n\choose k}\}_{k=0}^n$ satisfies the inequality (\ref{2}). \qed

Consequently, we have found a sufficient condition for both the $2$-log-concavity and the higher order Tur{\'a}n inequalities of the sequence $\{{n\choose k}\}_{k=0}^n$.

Recall that there is a simple criterion on a nonnegative sequence $\{a_k\}_{k=0}^{\infty}$ that guarantees infinite log-concavity. Namely
$$a_k^2\geq ra_{k-1}a_{k+1},$$
where $r\geq(3+\sqrt{5})/2$, for all $k\geq 1$.

For the inequality (\ref{2}), we have the following result.

\begin{thm}
 The positive sequence $\{a_k\}_{k=0}^{\infty}$ satisfies the inequality (\ref{2}), if
  \begin{equation}\label{inf}
  a_k^2\geq ra_{k-1}a_{k+1},
  \end{equation}
where $r\geq 2$, for all $k\geq 1$.
\end{thm}
\pf Applying the inequality (\ref{inf}), we have
\begin{equation}
\frac{a_k}{a_{k-1}}\geq 2\frac{a_{k+1}}{a_k}.
\end{equation}
Since $a_k\geq 0$ for $k\geq 0$, we easily get
\begin{equation}
\frac{a_{k+2}}{a_{k+1}}+\frac{a_k}{a_{k-1}}\geq 2\frac{a_{k+1}}{k},
\end{equation}
which means the sequence $\{\frac{a_{k+1}}{a_k}\}_{k\geq 0}$	is convex.

Similarly, we can prove that the sequence $\{\frac{a_k}{a_{k+1}}\}_{k\geq 0}$	is also convex. \qed

Notice that in the proof of Theorem \ref{2l}, we did not apply the Theorem \ref{mainthm}. In the last part of this section, we will prove the following result.

\begin{thm}\label{3main3}
	For a log-concave positive sequence $\{a_k\}_{k\geq0}$, if it satisfies the following inequalities
	\begin{equation}\label{333}
	\frac{a_{k+2}}{a_k}
	\leq\frac{a_{k+1}^2-a_ka_{k+2}}{a_k^2-a_{k-1}a_{k+1}}
	\leq \frac{a_{k+1}}{a_{k-1}},
	\end{equation}
	for $k\geq 1$. Then 
	\begin{equation}\label{334}
	\frac{a_{k+1}^2}{a_k^2}(1-\sqrt{1-c_k})^2\leq \frac{a_{k+2}}{a_k}
	\leq\frac{a_{k+1}^2-a_ka_{k+2}}{a_k^2-a_{k-1}a_{k+1}}
	\leq
	\frac{a_{k+1}^2}{a_k^2}(1+\sqrt{1-c_k})^2
	\end{equation}
where $c_k=a_ka_{k+2}/a_{k+1}^2$,	for $k\geq 1$.
\end{thm} 

\pf Since \begin{equation}
a_{k+1}^2-a_ka_{k+2}=a_{k+1}^2(1-c_k)
\end{equation} 
and $\{a_k\}_{k\geq0}$ is a positive sequence, the inequality (\ref{333}) is equivalent to
\[
\frac{a_{k+2}}{a_k}
\leq\frac{a_{k+1}^2(1-c_k)}{a_k^2(1-c_{k-1})}
\leq \frac{a_{k+1}}{a_{k-1}},
\]
i.e.
\begin{equation}\label{335}
c_k\leq \frac{1-c_k}{1-c_{k-1}}\leq \frac{1}{c_{k-1}}.
\end{equation}
And the inequality (\ref{334}) is equivalent to 
\begin{equation}\label{336}
(1-\sqrt{1-c_k})^2\leq c_k\leq \frac{1-c_k}{1-c_{k-1}}\leq (1+\sqrt{1-c_k})^2.
\end{equation}
We aim to prove the inequality (\ref{336}).

First we will prove that
\begin{equation}
(1-\sqrt{1-c_k})^2\leq c_k.
\end{equation}
 Since $\{a_k\}_{k\geq0}$ is a log-concave positive sequence, we have
$0\leq c_k\leq 1$.
Hence
\begin{equation}
0\leq \sqrt{1-c_k}\leq 1.
\end{equation}
Multiplying $\sqrt{1-c_k}$, we get
\begin{equation}
0\leq 1-c_k\leq \sqrt{1-c_k}.
\end{equation}
Note that $c_k=1-(1-c_k)$ and we obtain
\begin{equation}
c_k\geq 1-\sqrt{1-c_k}\geq (1-\sqrt{1-c_k})^2.
\end{equation}

Now we proceed to prove that 
\begin{equation}\label{339}
\frac{1-c_k}{1-c_{k-1}}\leq (1+\sqrt{1-c_k})^2.
\end{equation}
For $0\leq c_k\leq \frac{\sqrt{5}-1}{2}$, we have $1-c_k\geq \frac{3-\sqrt{5}}{2}$. Applying the right inequality of (\ref{335}), we get
\begin{equation}
\frac{3-\sqrt{5}}{2(1-c_{k-1})}\leq \frac{1}{c_{k-1}}.
\end{equation}
Hence $c_{k-1}\leq \frac{2}{5-\sqrt{5}}$. And $1-c_{k-1}\geq \frac{3-\sqrt{5}}{5-\sqrt{5}}$. It follows that
\begin{equation}
\frac{1-c_k}{1-c_{k-1}}\leq \frac{5-\sqrt{5}}{3-\sqrt{5}}(1-c_k).
\end{equation}
Therefore, it is sufficient to prove
\begin{equation}\label{337}
\frac{5-\sqrt{5}}{3-\sqrt{5}}(1-c_k)\leq (1+\sqrt{1-c_k})^2.
\end{equation}	
Set $t=\sqrt{1-c_k}$. Then $\frac{\sqrt{5}-1}{2}\leq t\leq 1$, $c_k=1-t^2$, and inequality (\ref{337}) is equivalent to 
\begin{equation}\label{338}
\frac{5-\sqrt{5}}{3-\sqrt{5}}t^2\leq (1+t)^2.
\end{equation}
It is easy to verify that the inequality (\ref{338}) holds for 
\begin{equation}
\frac{1}{2}(3-\sqrt{5}-\sqrt{20-8\sqrt{5}})\leq t\leq \frac{1}{2}(3-\sqrt{5}+\sqrt{20-8\sqrt{5}}).
\end{equation}
And we can verify that \begin{equation}\frac{1}{2}(3-\sqrt{5}-\sqrt{20-8\sqrt{5}})\leq\frac{\sqrt{5}-1}{2}\leq 1\leq \frac{1}{2}(3-\sqrt{5}+\sqrt{20-8\sqrt{5}}). 
\end{equation}
Consequently, we have finished the proof for the inequality (\ref{339}) for $0\leq c_k\leq \frac{\sqrt{5}-1}{2}$.

For $\frac{\sqrt{5}-1}{2}< c_k\leq 1$, actually we could prove that
\begin{equation}
\label{340}
\frac{1-c_k}{1-c_{k-1}}\leq\frac{1}{c_{k-1}}\leq (1+\sqrt{1-c_k})^2.
\end{equation}
Apply the left inequality of (\ref{335}), multiply $\frac{1-c_{k-1}}{c_k}$, and we get
\begin{equation}
1-c_{k-1}\leq \frac{1-c_k}{c_k}.
\end{equation}
It follows that
\begin{equation}
c_{k-1}\geq \frac{2c_k-1}{c_k},
\end{equation}
Thus
\begin{equation}
\frac{1}{c_{k-1}}\leq \frac{c_k}{2c_k-1}.
\end{equation}
We aim to prove that for $\frac{\sqrt{5}-1}{2}< c_k\leq 1$,
\begin{equation}\label{341}
\frac{c_k}{2c_k-1}\leq (1+\sqrt{1-c_k})^2.
\end{equation}
Set $t=\sqrt{1-c_k}$. Then $0\leq t\leq \frac{\sqrt{5}-1}{2}$, $c_k=1-t^2$, and inequality (\ref{341}) is equivalent to
\begin{equation}
\frac{1-t^2}{1-2t^2}\leq (1+t)^2.
\end{equation}
Multiplying $\frac{1-2t^2}{1+t}$ with both sides, we get
\[
1-t\leq (1+t)(1-2t^2),
\]
i.e.
\begin{equation}\label{342}
t(t^2+t-1)\leq 0.
\end{equation}
Obviously, the inequality (\ref{342}) holds for $t\leq -\frac{1+\sqrt{5}}{2}$ or $0\leq t\leq \frac{\sqrt{5}-1}{2}$.
Hence we complete the proof.\qed

\noindent{\it Remarks.}
  In fact, based on the Theorem \ref{3main3}, the Theorem \ref{2l} is a corollary of Theorem \ref{mainthm}. And in the proof above, for $0\leq c_k\leq \frac{\sqrt{5}-1}{2}$, we could not determine whether $\frac{1}{c_{k-1}}\leq (1+\sqrt{1-c_k})^2$ is true. We ask for an answer to this question.

\section{Further Work}
In this section, we want to discuss a problem we will concern in the future work.

For the partition function $p(n)$, Hou and Zhang \cite{hou} proved that $p(n)$ is $2$-log-concave for $n\geq 221$. The fact inspires us to consider the problem whether we can find a sufficient condition, similar to the inequality (\ref{2}), for both the $2$-log-concavity and the higher order Tur{\'a}n inequalities for $p(n)$ for $n\geq 221$.

Actually, Chen, Wang and Xie \cite{cwx} proved that $\{p(n+1)/p(n)\}_{n\geq 116}$ is log-convex. Hence $\{p(n+1)/p(n)\}_{n\geq 116}$ is convex.  Applying Lemma \ref{lem1}, then we get that 
\begin{equation}
\frac{p(k+1)^2-p(k)p(k+2)}{p(k)^2-p(k-1)p(k+1)}
\leq \frac{p(k+1)}{p(k-1)}
\end{equation}
holds for $k\geq 116$. However,  it seems that we can not find an integer $N\geq 0$ to make sure that the inequality 
\[
\frac{a_{k+1}^2-a_{k}a_{k+2}}{a_k^2-a_{k-1}a_{k+1}}
\geq \frac{a_{k+2}}{a_k}
\]
holds for $p(k)$ for $k\geq N$.

In deed,  set 
\begin{equation}\label{fn}
f(n)=\frac{p(n+1)^2-p(n)p(n+2)}{p(n)^2-p(n-1)p(n+1)},
\end{equation}
and
\begin{equation}\label{gkn}
g_k(n)=\frac{p(n+k+2)}{p(n+k)},
\end{equation}
then we can verify that, for $224\leq n\leq 225$,
\[
g_{20}(n)\leq f(n)\leq g_{20}(n-1),
\]
 for $244\leq n\leq 261$,
\[
g_{21}(n)\leq f(n)\leq g_{21}(n-1),
\]
for $268\leq n\leq 291$,
\[
g_{22}(n)\leq f(n)\leq g_{22}(n-1),
\]
for $296\leq n\leq 323$,
\[
g_{23}(n)\leq f(n)\leq g_{23}(n-1),
\]
for $326\leq n\leq 355$,
\[
g_{24}(n)\leq f(n)\leq g_{24}(n-1),
\]
for $356\leq n\leq 389$,
\[
g_{25}(n)\leq f(n)\leq g_{25}(n-1),
\]
and for $390\leq n\leq 425$,
\[
g_{26}(n)\leq f(n)\leq g_{26}(n-1).
\]

Based on the verification above, we guess for any integer $n\geq 326$, we can find the $k=k(n)$ to satisfies the inequality
\begin{equation}
\frac{p(n+k+2)}{p(n+k)}\leq \frac{p(n+1)^2-p(n)p(n+2)}{p(n)^2-p(n-1)p(n+1)}\leq \frac{p(n+k+1)}{p(n+k-1)}.
\end{equation}

For $k\geq 3$, we can easily prove that
\begin{equation}
\frac{p(n+k+1)}{p(n+k-1)}\leq \frac{p(n+1)^2}{p(n)^2}\left(1-\sqrt{1+\frac{p(n)p(n+2)}{p(n+1)^2}}\right)^2.
\end{equation}
Then if we can prove that
\begin{equation}
\frac{p(n+k+2)}{p(n+k)}\geq \frac{p(n+1)^2}{p(n)^2}\left(1-\sqrt{1-\frac{p(n)p(n+2)}{p(n+1)^2}}\right)^2,
\end{equation}
we will find the sufficient condition for both the $2$-log-concavity and the higher order Tur{\'a}n inequalities for $p(n)$ for $n\geq 326$.


\begin{thebibliography}{99} \small

\bibitem{boa} {\sc R. Boas}, {\em Entire Function}, Academic Press , New York, 1954.


\bibitem{bor} {\sc G. Boros and V. H. Moll}, {\em Irresistible Integrals}, Cambridge University Press, Cambridge, 2004.

\bibitem{bra} {\sc P. Br{\"a}nd{\'e}n}, {\em Iterated sequences and the geometry of zeros}, J. Reine Angew. Math.,
658(2011):115--131.


\bibitem{che} {\sc W.Y.C. Chen}, {\em The spt-Function of Andrews}, Surveys in Combinatorics 2017, A. Claesson, M. Dukes, S. Kitaev, D. Manlove and K. Meeks, Eds., Cambridge Univ. Press, Cambridge, 2017.

\bibitem{cdy} {\sc W.Y.C. Chen, D.Q.J. Dou, and A.L.B. Yang}, {\em Br{\"a}nd{\'e}n's conjectures on
the Boros-Moll polynomials}, Int. Math. Res. Not. IMRN, 20(2013):4819--4828.

\bibitem{cjw} {\sc W.Y.C. Chen, D.X.Q. Jia, L.X.W Wang}, {\em Higher Order Tur{\'a}n Inequalities for the Partition Function}, Trans. Amer. Math. Soc. 372 (2019) 2143--2165.

\bibitem{cwx} {\sc W.Y.C. Chen, L.X.W. Wang and G.Y.B. Xie}, {\em Finite differences of the logarithm of the partition function}, Math. Comp. 85 (298) (2016) 825--847.


\bibitem{cra} {\sc T. Craven and G. Csordas}, {\em Iterated Laguerre and Tur{\'a}n inequalities}, JIPAM. J.
Inequal. Pure Appl. Math., 3(3)(2002):Article 39, 14 pp. (electronic).

\bibitem{ccv} {\sc G. Csordas, T.S. Norfolk and R.S. Varga}, {\em The Riemann hypothesis and
theTur{\'a}n inequalities}, Trans. Amer.Math. Soc. 296 (2) (1986) 521--541.



\bibitem{dim1} {\sc D.K. Dimitrov}, {\em Higher order Tur{\'a}n inequalities}, Proc. Amer. Math.
Soc. 126 (7) (1998) 2033--2037.

\bibitem{dim2} {\sc D.K. Dimitrov and F.R. Lucas}, {\em Higher order Tur{\'a}n inequalities for the
Riemann $\xi$-function}, Proc. Amer. Math. Soc. 139 (3) (2011) 1013--1022.

\bibitem{fis} {\sc S. Fisk}, {\em Questions about determinants and polynomials}, arXiv: 0808.1850.

\bibitem{hon} {\sc R. Honsberger},  {\em More Mathematical Morsels}, Washington, DC: Math. Assoc. Amer., 237--239, 1991.

\bibitem{hou} {\sc Q-H. Hou, Z-R. Zhang},  {\em r-log-concavity of Partition Functions
	 }, Ramanujan J (2019) 48: 117. https://doi.org/10.1007/s11139-017-9975-5.

\bibitem{kau} {\sc M. Kauers and P. Paule}, {\em A computer proof of Moll's log-concavity conjecture},
Proc. Amer. Math. Soc. 135 (2007), 3847--3856.

\bibitem{lev} {\sc B.Ja. Levin}, {\em Distribution of Zeros of Entire Functions}, revised ed.,
Translations of Mathematical Monographs, vol.5, American Mathematical
Society, Providence, R. I., 1980.

\bibitem{mar} {\sc J. Ma{\u r}{\'i}k}, {\em On polynomials with all real zeros}, {\u C}asopis P{\u e}st. Mat. 89(1964), 5--9, MR 31:4782.


\bibitem{mcn} {\sc P.R.W. McNamara and B.E. Sagan}. {\em Infinite log-concavity: developments and
conjectures}, Adv. in Appl. Math., 44(1)(2010):1--15.

\bibitem{mil} {\sc G.V. Milovanovi{\' c}, D.S. Mitrinovi{\' c} and Th. M. Rassias}, {\em Topics in polynomials: extremal problems, inequalities, zeros}, World Scientific, Singapore, 1994. MR 95m:30009.

\bibitem{pol1} {\sc G. P{\'o}lya}, {\em  {\"U}ber die algebraisch-funktionentheoretischen Untersuchungen vonj. L. W. V. Jensen}, Kgl.
Danske Vid. Sel. Math.-Fys. Medd. 7 (1927), 3--33.

\bibitem{pol} {\sc G. P{\'o}lya and J. Schur}, {\em {\"U}ber zwei Arten von Faktorenfolgen in der
Theorie der algebraischen Gleichungen}, J. Reine Angew. Math. 144
(1914) 89--113.

\bibitem{rah} {\sc Q.I. Rahman and G. Schmeieer}, {\em Analytic Theroy of Polynomials}, Oxford
University Press, Oxford, 2002.

\bibitem{sta} {\sc R.P. Stanley}, {\em Personal communication}, May 2008.



\end{thebibliography}
\end{document}